\documentstyle{article}

\newcommand{\R}{{\bf R}}

\newcommand{\Z}{{\bf Z}}
\newcommand{\C}{{\bf C}}

\newcommand{\Sp}{{\rm Sp}}
\newcommand{\SU}{{\rm SU}}
\newcommand{\p}{{\partial}}

\newcommand{\Om}{{\Omega}}
\newcommand{\om}{{\omega}}
\newcommand{\eps}{{\varepsilon}}

\newcommand{\Ll}{{\cal L}}

\newcommand{\Nn}{{\cal N}}

\newcommand{\inner}[2]{\langle #1, #2\rangle}
\newcommand{\IMP}{\Longrightarrow}

\newcommand{\half}{{\scriptstyle\frac{1}{2}}}
\newcommand{\grad}{{\rm grad }}
\newcommand{\Si}{{\Sigma}}
\newcommand{\SS}{{\smallskip}}
\newcommand{\MS}{{\medskip}}
\newcommand{\NI}{{\noindent}}
\newcommand{\proof}[1]{\noindent{\bf Proof#1:\  }}

\newcommand{\QED}{\hfill$\Box$\medskip}

\begin{document}

\NI
{\bf\large
Introduction to Symplectic Topology: corrigenda}

Several readers have pointed out to us various small errors and typos
in this book. 
All are minor except for an error in the statement of
Theorem 3.17 on p 94 spotted by David Theret.  We 
thank him as well as everyone else who told us of these errors. 

The first part of this note is
a list of short corrections.  The second part contains
some longer revised passages.

\begin{flushright}
Dusa McDuff and Dietmar Salamon.  Sept 1997
\end{flushright}
\MS

\SS\SS
\NI
{\bf\large A list of Short Corrections}

\NI
p 10 line 10:\quad
``...such as a pendulum or top...''
\SS

\NI
p 17 line 7:\quad   $J_0 X_{H_t} = \nabla H_t$
\SS

\NI
p 22 line 20:\quad
``...their solutions.''\quad
line 28:\quad
``...function of the variables...''
\SS

\NI
p. 24 formula (1.19):  for consistency with later formulas  
there should be a $-$ sign in this equation
\SS

\NI
p. 39, line 5 of Proof:\quad  for all $v\in W$  
\SS

\NI
p. 43 
On line 14 the text should read:\quad  
``To prove this, we choose a positive definite and symmetric matrix
$P\in\Sp(2n)$ such that....''
and on line 21/22:\quad  
``...choosing $P$ is equivalent to choosing a $G$-invariant 
inner product on $\R^{2n}$ that is compatible with
$\om_0$ in the sense that it has the form 
$\om_0(\cdot, J\cdot)$ for some $\om_0$-compatible 
almost complex structure $J$.''
\SS

\NI
p 44 line 1/2:\quad
``...sends a matrix $U\in\SU(n)$ to...''
\SS

\NI
p. 47 line 18:\quad  $ \mu(\Psi) =  \frac 12 \sum_t
{\rm sign}\,\Gamma(\Psi,t)$ 
\SS

\NI
p. 50 line -3:\quad  delete repetition of ``intersection"
\SS

\NI
p. 51 line 10:\quad $ -\langle\dot X(t)u,Y(t)u\rangle$
\SS

\NI
p. 53 lines 5 and 6 of Proof:\quad replace 
$\om_0(\Psi^{\rm T}v,\Psi^{\rm T}w)$ by 
$\om_0(\Psi^{\rm T}u,\Psi^{\rm T}v)$ twice
\SS

\NI
p. 55 line 10 of Proof:\quad $A\bar z_j=-i\alpha_j z_j$
\SS

\NI
p. 60 line  5,6:\quad  replace $E^\om$ by the orthogonal complement 
$E^\perp$
\SS

\NI
p. 65 line 10:\quad  delete
repetition of ``that"
\SS

\NI
p. 71 Exercise 2.66:\quad  
It would be more clear to say: ``there are 
precisely two" instead of ``one nontrivial"
\SS

\NI
p. 77 last line of Exercise 2.72:\quad  $c_1(\nu_{\C P^1})=1$.
\SS

\NI
p. 79, line -10:\quad
``Nondegeneracy...''
\SS

\NI
p. 82  line -1 of proof: \quad
$\iota({\psi_t}^*X)\om={\psi_t}^*(\iota(X)\om)$\SS
\SS

\NI
p. 86, line 2:\quad
``By Exercise~2.15, every such...''
\SS

\NI
p. 89, line 2:\quad
``...any vector $v^*\in T_q^*L$ can...''
\SS

\NI
p. 94 Statement of Theorem 3.17:\quad  Delete the last sentence.  
As David Theret pointed out, it is easy to find a counterexample to this
statement unless one requires that the class $[\om_t] - [\om_0]\in
H^2(M,Q;\R)$ is constant.  \SS

\NI
p. 96 line -4 of proof: replace $N$ by $\Nn$.\qquad
last 3 lines of Proof:\quad Name the diffeomorphisms $\chi_t$ instead
of $\phi_t$.   Then $\chi_t^*\om_t = \om_0 = \om$ and the desired extension
is   $\rho_t\circ\chi_t$.\SS

\NI
p. 100 line 2:\quad   $\phi: \Nn(L_0)\to V$\SS

\NI
p. 103 line -16:\quad
``...that $d\alpha$ restricts to...''
\qquad line -2:\quad  
Replace Corollary 2.4 by Corollary 2.5
\SS

\NI
p. 111 line -3:\quad
``...symplectization of $Q$.''
\SS

\NI
p. 113 line -2 of Proof:\quad 
$\psi^*\om = e^\theta(d\alpha-\alpha\wedge d\theta)$\SS

\NI
p. 115 line 14:\quad
``...$f:\C^{n+1}-\C^n\times\{-1\}\to\C^{n+1}-\C^n\times\{-i\}$...''
\SS

\NI
p. 117 line -6/-5:\quad
``...metric $g(u,v)=\inner{u}{v}$...''
\SS

\NI
p. 123 line 8:\quad
boldface ``{\bf (iv)}''.
\SS

\NI
p. 129 line 2: \quad  $\lambda: M\to \R$ 
\SS

\NI
p 167 line -5:\quad  
replace ``$p_1\sim p_2$'' by ``$p_0\sim p_1$''. 
\SS

\NI
p 169 line -3:\quad  $T_p{\cal O}(p)\subset (T_p(\mu^{-1}(0)))^\om.$
\SS

\NI
p. 172 line -10:\quad  delete the repetition of ``of''
\SS

\NI
p. 215 line 3 in Lemma~6.31:\quad  
bracket missing in ``$\inner{c_1(\nu_\Sigma)}{[\Sigma]}$''.
\SS

\NI
p. 219 Exercise 6.38:\quad the curve $C_3$ should be $\{z_1 = {z_2}^3\}$.
\SS

\NI
p. 221 Lemma 6:40:\quad  
    ``...$L(\delta)-L(0)$ is symplectomorphic to 
       the spherical shell $B(\lambda+\delta')-B(\lambda)$ for
      $\delta'=\sqrt{\lambda^2+\delta^2}$.''
\SS

\NI
p. 240 last line in definition of $\tilde\tau_t(x)$:\quad
Replace ``$c-\eps\le f(x)$'' by ``$f(x)\le c-\eps$''.
\SS

\NI
p. 253 lines 8,9:\quad 
``...that, up to diffeomorphism,  there is...''
\SS

\NI
p 265 The proof of Lemma~8.2 must be revised (see below). 
\SS

\NI
p 273 line -3:\quad
Replace ``$\R^{2n}$'' by ``$\R^{2}$''.
\SS

\NI
p 274 line 1/2:\quad
Replace ``$\R^{n}$'' by ``$\R$'' and
``$\R^{2n}$'' by ``$\R^{2}$''.
\SS

\NI p 303 line 8:\quad  
$...F_t = \int_0^t ...$ (not $...F_t=\int_0^1...$);
\newline
line -4:\quad
``$\Lambda={\rm graph}(dS)$.''
\SS

\NI
p 304 line 2:\quad  
``...the above action function...''
\SS

\NI
p 306 line 12:\quad  
``...where $A=A^T=\p_x\p_x\Phi\in\R^{n\times n},...$''
\newline
line -8:\quad
replace ``$Ax'-B\xi'$'' by ``$Ax'+B\xi'$''.
\SS

\NI
p 337 line 14:\quad  
replace ``$\cup a_k$'' by ``$\cup a_N$''.
\SS

\NI
p 340 line -11:\quad  
``...every $t$.  A compact invariant set...''
\SS

\NI
p. 360  line -3:\quad  $\psi$ is a symplectic embedding not a
symplectomorphism\SS

\NI
p. 361 line 16:\quad replace $\overline{c}_G$ by $\overline{w}_G$\SS

\NI
p 362 line 1:\quad  $c(E) = \pi r_1^2 = w_L(E)$\SS

\NI
p 364 line 5 of Proof:  the $\psi_t$ are diffeomorphisms not
symplectomorphisms\SS

\NI
p. 367 lines -7, -3:\quad
$\Ll(\{\phi_t\})$\SS

\NI
p 374 line -12:\quad  
$f$ should be a smooth embedding rather than a
 diffeomorphism
\SS

\NI
p. 376 line 6:\quad   the restriction of $\psi_H$ to $Z_{2c}$ 
is called  $\Psi_H$\qquad

line -1:\quad  $\pi R^2 = c + e + \epsilon$
\SS

\NI
p. 377 line 3,4 of Exercise 12.22:\quad 
The text should read:\quad
``... symplectic embedding of the ball $B^{2n+2}(r)$ into
$B^2(R)\times M$ where $\pi R^2 = e + c/2+\epsilon$.
Using the...''
\SS

\NI
p 378 line 3\quad: the supremum and infimum 
should be taken over $x\in\R^{2n}$
\SS

\NI
p. 384 line 2 of Proof:\quad  $\Ll_X\om_0 = \om_0$
\SS

\NI
p 390 line -11:\quad the conditions (I), (II), (III).
\SS

\NI
p 398 line 1:\quad 
$V_\tau(x_{j+1}+s\xi_{j+1},y_j+s\eta_j) $
\SS

\MS\MS
\NI
{\bf \large Longer revised passages}
\MS

\NI page 72: Remark 2.68~(ii)

Several students have pointed out that the sentence ``The axioms
imply that this integer depends only on the homology class of $f$"
is hard to substantiate. 
Change this and the rest of Remark 2.68~(ii) to:
\MS

``We will see in Exercise~2.75
this integer depends only on the homology class
of $f$. Thus the first Chern number generalizes to
a homomorphism $H_2(M,\Z)\to\Z$.  This gives rise 
a cohomology class $c_1(E)\in H^2(M,\Z)/{\rm torsion}$.
There is in fact a natural choice of a lift of 
this class to $H^2(M,\Z)$, also denoted by $c_1(E)$, 
which is called the {\bf first Chern class}.
We shall not discuss this lift in detail, 
but only remark that in the case of a line bundle 
$L\to M$ the class $c_1(L)\in H^2(M,\Z)$ 
is Poincar\'e dual to the homology class
determined by the zero set of a generic section.''
\MS

Then replace the current exercises 2.75 and 2.76 
by the following new version:
\MS

\NI {\bf Exercise 2.75}
{\bf (i)} Prove that every symplectic vector 
bundle $E$ over a Riemann surface $\Si$ decomposes
as a direct sum of $2$-dimensional 
symplectic vector bundles.
{\bf Hint:}  Show that any such vector bundle
of rank $> 2$ has a nonvanishing section.

\NI{\bf (ii)} Suppose that $\Si$ is oriented
and that the bundle $E$
above extends
over a compact oriented $3$-manifold
$X$ with  boundary $\p X=\Sigma$.  
Prove that the restriction $E|\Sigma$
has Chern class zero.
{\bf Hint:}  Use (i) above and
look at a section $s$ as in Theorem~2.71.

\NI{\bf (iii)} Use (i) and (ii) above to substantiate
the claim made in Remark~2.68
above that the Chern class $c_1(f^*E)$ 
depends only on the homology
class of $f$.  Here the main problem is that
when $f_*([\Si])$ is null-homologous the 
$3$-chain that bounds it 
need not be representable by a $3$-manifold.
However its singularities can be assumed to have
codimension $2$ and so the proof of (ii) goes 
through.  
\QED

\MS

\MS\MS\MS

This is a new exercise that should go at the very
end of Chapter 2.
\MS

\NI
{\bf Exercise}
Prove that every symplectic vector bundle 
$E\to\Sigma$ over a Riemann surface $\Sigma$
which admits a Lagrangian 
subbundle can be symplectically trivialized.
{\bf Hint:} Use the proof of Theorem~2.67
to show that $c_1(E)=0$.
\QED

\MS\MS

This is a revised version of Lemma 8.2.
\MS

\NI
{\bf Lemma 8.2 } {\it
The Poincar\'e section $\Sigma\cap U$ is a symplectic
submanifold of $M$ and the Poincar\'e section map 
$\psi:\Sigma\cap U\to\Sigma$ is a symplectomorphism.}
\MS

\proof{}  The hypersurface $\Sigma$ is of dimension $2n-2$ and the tangent
space at $p$ is
$$
    T_p\Sigma = \{v\in T_pM\,|\,dG(p)v = dH(p)v = 0\}.
$$
The condition $\{G,H\} = \om(X_G,X_H)\neq 0$ shows that the
$2$-dimensional subspace spanned by $X_G(p)$ and $X_H(p)$ is a
complement of $T_p\Sigma$.  Now let $v\in T_p\Sigma$ and suppose that
$\om(v,w) = 0$ for all  $w\in T_p\Sigma$.  Then $\om(X_H(p),v) =
dH(p)v = 0$ and  $\om(X_G(p),v) = dG(p)v = 0$ and hence $v = 0$.  
Thus the $2$-form $\om$ is non\-degenerate on the subspace
$T_p\Sigma\subset\R^{2n}$.

To prove that $\psi$ is a symplectomorphism we consider the $2$-form
$$
     \om_H = \om + dH\wedge dt 
$$
on $\R\times M$. This is the {\bf differential form of Cartan}. 
\index{Cartan!differential form of} 
It has a $1$-dimensional kernel consisting of 
those pairs $(\theta,v)\in\R\times T_pM$ which satisfy 
$$
     v = \theta X_H(p).
$$
Now let $D\subset\C$ denote the unit disc in the complex plane and let    
$u:D\to\Sigma$ be a $2$-dimensional surface in $\Sigma$.
We must prove that
$$
     \int_{D} u^*\psi^*\om = \int_D u^*\om.
$$
To see this consider the manifold with corners
$$
      \Om=\left\{(t,z)\,|\,z\in D,\,0\le t\le \tau(u(z))\right\}
$$
and define $v:\Om\to\R\times M$ by 
$$
     v(t,z) = \left(t,\phi^t(u(z))\right).
$$
Denote $v_0(z)=v(0,z)$ and $v_1(z)=v(\tau(u(z)),z)$.
Then ${v_0}^*\om_H=u^*\om$ and ${v_1}^*\om_H=u^*\psi^*\om$. Moreover, 
the tangent plane to the surface $v(\R\times\p D)$ contains the kernel of
$\om_H$.  Hence the $2$-form $v^*\om_H$ vanishes on the
surface   $\R\times\p D$.   Since $\om_H$ is closed it follows
from Stokes' theorem that 
$$
     0 = \int_{\Om} v^*d\om_H 
       = \int_{\p\Om}\om_H
       = \int_{D} u^*\psi^*\om - \int_D u^*\om.
$$
Hence $\psi$ is a symplectomorphism. \QED
\MS\MS

This is a revised version of Lemma 12.37 and Exercise 12.38.
\MS

\NI
{\bf Lemma 12.37 } {\it
Let $H$ be any Hamiltonian which equals 
$$
    H_\infty(z)  = (\pi + \eps)|z_1|^2+\half {\pi}|z_r|^2.
$$ 
for large $|z|$. Then the functional
$\Phi_H^\tau:\R^{2nN}\to\R$ satisfies the Palais--Smale condition.}

\MS

\proof{}  The Palais--Smale condition asserts 
that for every sequence ${\bf z}^\nu$ in $\R^{2nN}$
$$
       \left\|\grad\,\Phi_H^\tau({\bf z}_\nu)\right\|_\tau\to 0
       \qquad\IMP\qquad
       \sup_\nu\left\|{\bf z}_\nu\right\|_\tau<\infty.
$$
Suppose otherwise that
$\left\|{\bf z}_\nu\right\|_\tau\to\infty$.  Then, since 
$\grad\,\Phi_H^\tau({\bf z}_\nu)$ converges to zero, we claim that
all components $z_j^\nu$ of ${\bf z}_\nu$ must diverge to infinity.
Clearly, this will follow if we prove the inequality
$$
    \min_j|z_j^\nu|\ge \frac{1}{c}\max_j|z_j^\nu| - 1
$$
for some constant $c\ge 1$ which is independent of $\nu$.
A proof of this is sketched in Exercise~12.38 below. 
Hence we may assume that the $z_j^\nu$ all 
lie in a region in which  $H(z)=H_\infty (z)$.
Now consider the sequence
$$
       {\bf w}_\nu 
       = \frac{{\bf z}_\nu}
         {\left\|{\bf z}_\nu\right\|_\tau}.
$$
This sequence has a convergent subsequence, and it is easy to check
that the limit  has norm $1$ and is a critical point of
$\Phi_{H^\infty}^\tau$.  But, because the flow of $H_\infty$ has no
non\-constant periodic orbits of period $1$, the fixed point $0$ is the
only critical point of this functional.  This contradiction proves the
lemma. \QED
\MS\MS

\noindent

\NI {\bf Exercise 12.38}
This exercise fills in a missing detail in the proof of the
above Lemma.  
Given a vector ${\bf z}$ in $\R^{2nN}$ with components 
$z_j=(x_j,y_j)\in\R^{2n}$ denote by $\zeta_j=(\xi_j,\eta_j)\in\R^{2n}$  
the $j$th component  of $\grad\,\Phi_H^\tau({\bf z})$.  
Then $x_{j+1}$ is the unique solution of the equation 
$$
    x_{j+1} = F_j(x_{j+1},\eta_j)
$$
where 
$$
       F_j(x_{j+1},\eta_j) 
       = x_j 
         + \tau\frac{\p V_\tau}{\p y}(x_{j+1},y_j) 
         + \tau\eta_j. 
$$
Prove that for sufficiently small $\tau$ and any $\eta_j$ the map
$x_{j+1}\mapsto F_j(x_{j+1},\eta_j)$ is a contraction 
with Lipschitz constant 
$
     \alpha 
     = \tau \sup_{x}\left|\p^2 V_\tau/\p x\p y(x,y_j)\right|
     < 1.
$
Deduce that 
$$
     |x_{j+1}-x_j| 
     \le \frac{\tau}{1-\alpha}
         \left|\frac{\p V_\tau}{\p y}(x_j,y_j)+\eta_j\right|.
$$
Use this and the inequality
$$
     |y_{j+1}-y_j| 
     \le \tau
         \left|\frac{\p V_\tau}{\p x}(x_{j+1},y_j)+\xi_{j+1}\right|
$$
to conclude that, if $\|\grad\,\Phi_H^\tau({\bf z})\|_\tau \le 1$
and $\tau$ is sufficiently small then
$$
     |z_{j+1}-z_j|\le \frac{1}{2}(|z_j|+1).
$$
This implies 
$$
     2^{-j}(|z_0|+1) \le |z_j|+1 \le 2^j(|z_0|+1)
$$
for $j=0,\dots,N$.  Hence, if 
${\bf z}_\nu\in\R^{2nN}$ is a sequence with 
$\|\grad\,\Phi_H^\tau({\bf z}_\nu)\|_\tau \le 1$
and $\tau$ sufficiently small such that 
$\left\|{\bf z}_\nu\right\|_\tau\to\infty$,
then $\left\|z_{\nu j}\right\|_{\R^{2n}}\to\infty$
for all $j$. 
\QED

\MS\MS

\NI
Dusa McDuff:  dusa@math.sunysb.edu

\NI
Dietmar Salamon: das@maths.warwick.ac.uk

\end{document}